\documentclass[journal]{IEEEtran}
\ifCLASSINFOpdf
\else
\fi

\usepackage{balance}
\usepackage{float}
\usepackage{caption}
\usepackage{subcaption}
\usepackage{tikz,pgfplots,pgfplotstable,booktabs}
\usepackage[official]{eurosym}
\usepackage{xcolor}
\usepackage{colortbl}
\usepackage{comment}
\usepackage{circuitikz}
\usepackage{amsmath,amsfonts,amsthm,bm}
\allowdisplaybreaks[1]

\pgfplotsset{
    discard if not/.style 2 args={
        x filter/.code={
            \edef\tempa{\thisrow{#1}}
            \edef\tempb{#2}
            \ifx\tempa\tempb
            \else
                
            \fi
        }}}

\begin{document}
%
\title{\textcolor{black}{Data-Driven Screening of Network Constraints for Unit Commitment}}
%
%
%

\author{S. Pineda, J. M. Morales and A. Jim\'enez-Cordero

\thanks{S. Pineda is with the Department
of Electrical Engineering, University of M\'alaga, M\'alaga, Spain. E-mail: spinedamorente@gmail.com.}
\thanks{J. M. Morales and A. Jim\'enez-Cordero are with the Department of Applied Mathematics, University of M\'alaga, M\'alaga, Spain. E-mails: juan.morales@uma.es; asuncionjc@uma.es.}
\thanks{This work was supported in part by the Spanish Ministry of Economy, Industry, and Competitiveness through project ENE2017-83775-P; and in part by the European
Research Council (ERC) under the EU Horizon 2020 research and innovation
program (grant agreement No. 755705). The authors thankfully acknowledge the computer resources, technical expertise and assistance provided by the SCBI (Supercomputing and Bioinformatics) center of the University of M\'alaga.}}

\maketitle

\begin{abstract}
The transmission-constrained unit commitment (TC-UC) problem is one of the most relevant problems solved by independent system operators for the daily operation of power systems. Given its computational complexity, this problem is usually not solved to global optimality for real-size power systems. In this paper, we propose a data-driven method that leverages historical information to screen out network constraints in the TC-UC problem. First, past data on demand and renewable generation throughout the network are used to \textit{learn} the congestion status of transmission lines. Then, we infer the lines that will not become congested for upcoming operating conditions based on such learning and disregard their capacity constraints.  This way, we formulate a reduced TC-UC problem that is easier to solve. Numerical results on a medium- and a large-size power system show that the proposed approach outperforms existing ones by significantly reducing the computational time while obtaining solutions that are equal or close to the one obtained with the original TC-UC problem. Furthermore, the purely data-driven method we propose can be seamlessly complemented with a constraint generation procedure to guarantee that the optimal solution to the original TC-UC problem is eventually recovered.
\end{abstract}


%
\IEEEpeerreviewmaketitle

\section{Introduction} \label{sec:introduction}

\IEEEPARstart{D}{espite} the liberalization of the electricity sector, the centralized transmission-constrained unit-commitment problem (TC-UC) is still a crucial optimization tool for market clearing and the operation of current power systems \cite{Zheng2015}. More specifically, the goal of TC-UC is to determine the commitment and dispatch of generating units to satisfy electricity demand at the minimum cost. The general mathematical formulation of the TC-UC problem can be stated as follows: 
\begin{subequations}
\begin{align}
\min_{\textbf{x}\in \mathbb{R}^n,\,\textbf{y}\in\{0,1\}^m} & \quad f(\textbf{x},\textbf{y}) \label{uc1_of}\\   
 & \quad g_i(\textbf{x},\textbf{y}) \leq 0, \quad \forall i \label{uc1_c1} \\
 & \quad h_j(\textbf{x}) \leq 0, \quad \forall j \label{uc1_c2} 
\end{align}\label{uc1}%
\end{subequations}
where continuous variables $\textbf{x}$ represent power dispatches, power flows through lines, unserved load, etc. and binary variables $\textbf{y}$ model the on/off status of the generating units. 

The objective function \eqref{uc1_of} minimizes the total generation costs, while technical constraints pertaining to generating units and the transmission network are modeled through \eqref{uc1_c1} and \eqref{uc1_c2}, respectively. If the objective function and all constraints are assumed linear, formulation \eqref{uc1} is a mixed-integer linear programming problem (MILP). Using a classical reduction from the knapsack problem, the authors of \cite{Bendotti2019} show that the TC-UC is NP-hard even if a single time period is considered. One possible approach to reduce the computational burden of the TC-UC problem is to define new efficient reformulations of the MILP problem, as done for instance in \cite{Carrion2006, Morales-Espana2013}. Unfortunately, the classical solution algorithms proposed in such works fail to provide solutions within a reasonable amount of time for real-size power systems~\cite{Sioshansi2008a}.

In a different line of thought, some authors have suggested the use of statistical learning algorithms to provide a quick proxy of commitment decisions \cite{Dalal2018}. Similarly, reference \cite{Ng2018} presents a learning method to obtain optimal solutions of the optimal power flow with high probability.
From a market perspective, these approaches lack transparency and thus, its use seems more appropriate for long-term planning.

A more suitable alternative to reduce the computational burden of unit-commitment problems for operating applications consists in removing superfluous technical constraints. Following this line of argument, the authors of \cite{Ardakani2013} introduce the concept of \textit{umbrella constraint} and propose a convex optimization problem to discover such constraints and reduce the complexity of the security-constrained unit-commitment problem. However, as stated by the authors, the solution times of the umbrella discovery problem are prohibitive for practical power system use. Inspired by the constraint generation approach \cite{Ben-Ameur2006}, a method that iteratively includes violated security constraints to the unit commitment problem is proposed in \cite{Fu2013}. Reference \cite{Madani2017} describes a cheap and parallelizable computational method to eliminate redundant constraints under the assumption that security constraints can be estimated as linear constraints. Finally, authors of \cite{Zhai2010} establish a necessary and sufficient condition for a security constraint to be inactive.

\begin{figure}
    \centering
    \begin{tikzpicture}[scale=0.8,every node/.style={scale=0.8}]
	\begin{axis}[xmin=0,xmax=5,ymin=0,ymax=5,xlabel=$x$,ylabel=$y$]
	\addplot[lightgray,ultra thick] coordinates {(1,3)(4,3)};
	\addplot[lightgray,ultra thick] coordinates {(2,2)(4,2)};
	\addplot[lightgray,ultra thick] coordinates {(3,1)(4,1)};
	\addplot[thick] coordinates {(0.5,3.5)(4,0)};
	\addplot[thick] coordinates {(4,0)(4,3.5)};
	\addplot[thick] coordinates {(4,3.5)(0.5,3.5)};
	\addplot[thick] coordinates {(0,4.5)(5,4.5)};
	\addplot[thick,dashed] coordinates {(4,3.5)(4,5)};
	\addplot[mark=o,thick] coordinates {(4,3)};
	\addplot[mark=o,thick] coordinates {(4,4)};
	\node[] at (axis cs: 4.2,3) {$A$};
	\node[] at (axis cs: 4.2,4) {$B$};
	\node[] at (axis cs: 4.3,1.5) {(2b)};
	\node[] at (axis cs: 1.7,1.7) {(2c)};
	\node[] at (axis cs: 2.5,4.7) {(2d)};
	\node[] at (axis cs: 2.5,3.7) {(2e)};
    \end{axis}
\end{tikzpicture}
    \caption{Type of constraints in a MILP}
    \label{fig:ilustration_mip}
\end{figure}
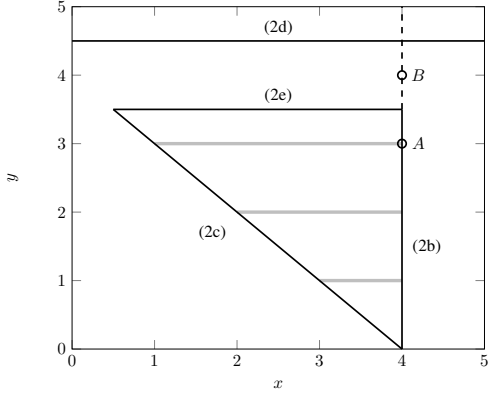

In this paper, instead, we focus on removing \emph{network} constraints in the TC-UC problem. As explained in \cite{Ostrowski2012}, network constraints \eqref{uc1_c2} significantly increase the time needed to solve the relaxed linear problems at each node of the branch-and-bound tree. Since a large amount of relaxed linear problems have to be solved to find the solution to \eqref{uc1}, removing some constraints \eqref{uc1_c2} in all relaxed linear problems may involve a significant time reduction. However, given the NP-hardness of mixed-integer linear problems, the extent of this reduction is not straightforward. Within the particular context of removing network constraints in the TC-UC problem, we highlight the more recent publication \cite{Roald2019}, which proposes a method to identify line flow constraints that do not become active for a large range of load variation so that the results are valid for long periods of time. To better illustrate our contributions with respect to the state-of-the-art, we consider the following illustrative MILP:
\begin{subequations}
\begin{align}
    \max_{x\in\mathbb{R},y\in\mathbb{Z}} & \quad x + y \label{ex_of}\\
    \text{s.t.} & \quad x \leq 4 \label{ex_c1}\\
    & \quad x+y \geq 4 \label{ex_c3}\\
    & \quad y \leq 4.5 \label{ex_c4} \\
    & \quad y \leq 3.5 \label{ex_c2}
\end{align} \label{ex}%
\end{subequations}
The four constraints of \eqref{ex} are depicted in Fig.~\ref{fig:ilustration_mip}, in which feasible points are represented by the solid light gray lines. It is just apparent that the optimal solution to problem \eqref{ex} is $A=(4,3)$. Let us now analyze the role of the four constraints:
\begin{itemize}
    \item[-] Constraint \eqref{ex_c1} holds with equality at the optimum. If this constraint is removed both the feasible region and the solution change. In fact, problem \eqref{ex} becomes unbounded. We refer to this type of constraint as \textit{active constraint}.
    \item[-] Constraint \eqref{ex_c3} does not hold with equality at the optimal point. If this constraint is removed, the feasible region changes but the optimal solution remains the same (for this particular objective function). We refer to this type of constraint as \textit{inactive constraint}.
    \item[-] Constraint \eqref{ex_c4} does not hold with equality at the optimal point either. Note, however, that this is a very particular \textit{inactive} constraint since, if removed, the feasible region is unaffected. This means that the optimal solution for \emph{any} given objective function remains the same when this constraint is eliminated. We refer to this type of constraint as \textit{redundant} constraint.
    \item[-] Constraint \eqref{ex_c2} is is not binding at the optimum either. However, if removed, the optimal solution changes from $A$ to $B=(4,4)$. Then, this is the only non-binding constraint that cannot be eliminated without affecting the optimal solution. We refer to this type of constraint as \textit{quasi-active constraint}.
\end{itemize}

In order to formulate a simplified problem with a lower computation burden than the original one, it is desirable to remove the \emph{inactive} \eqref{ex_c3} and \emph{redundant} \eqref{ex_c4} constraints. 
%

Therefore, the aim of this paper is to develop a novel data-driven method which efficiently identifies both \emph{redundant} and \emph{inactive}  network constraints of the TC-UC problem. Removing such constraints yields a simplified optimization problem that is faster to solve. Among the previously mentioned research works that propose methods to select a subset of network constraints that can be removed from \eqref{uc1}, references \cite{Fu2013, Zhai2010, Roald2019} are the ones closest to our work. On the one hand, the constraint generation method of \cite{Fu2013} removes both \textit{redundant} and \textit{inactive} constraints but requires to solve optimization problem \eqref{ex} more than once. On the other hand, the methods proposed in \cite{Zhai2010, Roald2019} aim at identifying \textit{redundant} constraints only. Consequently, these last two methods have three things in common. First, after removing these constraints, the two methods yield simplified optimization problems with the very same feasible region as the original problem. Second, the two methods are optimization-based and require the solution of additional optimization problems to identify \textit{redundant} constraints. This extra computational effort may turn these methods impractical for power system operation problems. Third, the impact of system congestion on the performance of those methods is not analyzed in either of the two references.

In this paper we propose a novel data-driven method to efficiently identify both \textit{redundant} and \textit{inactive} constraints, i.e., we take advantage of the information provided by the data to reduce the computational burden of the TC-UC problem. Therefore, the contributions of this paper are fourfold:
\begin{itemize}
    \item[-] Unlike \cite{Zhai2010,Roald2019}, the proposed method does not only identify \textit{redundant} constraints but also \textit{inactive} constraints such as \eqref{ex_c3}. In doing so, we remove a larger number of constraints from the original problem and consequently, reduce even further the computational time required to solve the simplified problem.
    \item[-] Unlike \cite{Fu2013,Zhai2010,Roald2019}, the core of the method we propose relies on learning from data instead of on solving a larger number of optimization problems. That is, we leverage available information of the power system to identify \textit{redundant} and \textit{inactive} constraints using a learning-based method that is computationally inexpensive.
    \item[-] The proposed purely data-driven approach can be seamlessly complemented, when necessary, with a constraint-generation procedure in the line of \cite{Fu2013}. By doing so, the obtained solution is guaranteed to be the same as that of the original problem, at the expense of a little increment in the computational time.
    \item[-] We compare the performance of the proposed approach with other existing methods in terms of computational burden and accuracy of a simplified network-constraint unit-commitment problem. For the sake of comparison, we investigate their performance for different levels of network congestion and system size.
\end{itemize}

The rest of this paper is organized as follows. Section \ref{sec:unit_commit} elaborates on the modeling assumptions and formulates the TC-UC problem. All the methods compared to determine the removed network constraints are explained in Section \ref{sec:methodology}. The comparison framework is described in Section \ref{sec:comparison}. Section \ref{sec:example} illustrates the different methods using a stylized example. The results of a medium- and large-size case studies are discussed in Section \ref{sec:num_results}. Finally, Section \ref{sec:conclusions} concludes the paper.

\section{Unit commitment problem} \label{sec:unit_commit}

In this section, we formulate a generic single-period TC-UC problem according to the following simplifying assumptions:
\begin{itemize}
    \item[-] \textit{Single-period}: Unit commitment is usually formulated as a multi-period problem that incorporates inter-temporal constraints such as ramping limits and minimum up and down times \cite{Carrion2006, Morales-Espana2013}. However, since this work focuses on investigating the impact of network reduction (\textit{spatial dimension}) on TC-UC results, we prefer to investigate such an impact alone by considering a single-period (\textit{temporal dimension}) TC-UC problem as in \cite{Roald2019}. 
    \item[-] \textit{DC power flow}: In order to keep the model linear, the power flows through the transmission network are computed using a DC approximation via power transfer distribution factors (PTDF). The PTDF of line $l$ with respect to node $n$ is denoted by $a_{ln}$. Besides, $\overline{f}_l$ represents the capacity of each transmission line. The number of buses and lines are denoted by $N_B$ and $N_L$, respectively.
    \item[-] \textit{Generation portfolio}: Each generating unit $g$ connected to node $b_g$ is characterized by a minimum and a maximum power output ($\underline{p}_g,\overline{p}_g$), a capacity factor ($\rho_g$) and a production cost ($c_g$). Thermal generating units ($\mathcal{G}^T$) have $c_g>0$, $\underline{p}_g>0$ and $\rho_g=1$, while renewable generating units ($\mathcal{G}^R$) have $c_g=0$, $\underline{p}_g=0$ and $0\leq \rho_g \leq 1$.
    \item[-] \textit{Known demand}: Electricity demand at each node $d_n$ is assumed to be known with certainty. 
    \item[-] \textit{No failures}: We assume that all units and lines are operational and thus security constraints are not considered.
\end{itemize}

Therefore, the single-period TC-UC  is formulated as the following MILP:
\begin{subequations}
\begin{align}
    \min_{p_{g},u_{g},q_{n},\epsilon_{n}}   & \quad \sum_g c_gp_{g} + L \sum_n |\epsilon_{n}|  \label{uc2_of}\\
     \text{s.t.} & \quad q_{n} + \epsilon_{n} = \sum_{g:b_g=n} p_{g} - d_{n}, \quad \forall n \label{uc2_net}\\
    & \quad \sum_n q_{n} = 0 \label{uc2_bal}\\
    & \quad u_{g} \underline{p}_g \leq p_{g} \leq u_{g} \rho_{g} \overline{p}_g, \quad \forall g \label{uc2_gen}\\
    & \quad -\overline{f}_l \leq \sum_n a_{ln} q_{n} \leq \overline{f}_l, \quad \forall l \label{uc2_cap}\\
    & \quad u_{g} \in \{0,1\}, \quad \forall g \label{uc2_bin}
\end{align} \label{uc2}%
\end{subequations}
Decision variables include the commitment of the generating units $u_g$, the power output dispatches $p_g$, the net power injections at each node $q_n$ and the slack variables $\epsilon_{n}$. Objective function \eqref{uc2_of} minimizes the total production cost to satisfy demand plus a penalty term of the slack variables. $L$ is a large enough constant such that $\epsilon_{n}=0, \forall n$ for a capacity-adequate power system. Equation \eqref{uc2_net} determines the net injection at each node, while constraint \eqref{uc2_bal} ensures power balance in the system. Constraints \eqref{uc2_gen} and \eqref{uc2_cap} impose limits on generation outputs and power flows using PTDF, respectively. Finally, binary variables are declared in \eqref{uc2_bin}.


\section{Methodology} \label{sec:methodology}

In this section, we present eight different methods to infer the subset of network constraints \eqref{uc2_cap} that can be removed from optimization problem \eqref{uc2}. For simplicity, the two inequality constraints in \eqref{uc2_cap} are either removed or kept simultaneously by all methods. The first four methods are simple benchmark methods whose results are not interesting \textit{per se} but serve us to measure the performance of the other approaches. The fifth, sixth and seventh methods are those proposed in references \cite{Fu2013}, \cite{Zhai2010} and \cite{Roald2019}, respectively. The last one is the data-driven method proposed in this paper.

Some methods rely on the availability of historical data from the power system under study. Historical data, also referred to as \textit{training} data set, is indexed by $t=1,\dots,N_T$. We assume we have access to historical data on i) hourly electricity demand at each bus $d_{n}$, ii) hourly capacity factors of all generating units, represented by $\rho_{g}$, and iii) the status of each transmission line $s_{l}$, which is equal to 1 if line $l$ was congested, and to 0, otherwise. The performance of the different methods are compared on a set of unseen operating conditions, also known as \textit{test} data set, indexed by $\hat{t}=1,\dots,\widehat{N_T}$.  

\subsection{Benchmark method (BN)} \label{sec:ben}

In the benchmark method, no network constraints are removed, i.e., the optimization problem is solved considering all the constraints. Hence, if used for real-size systems, the computational burden of the benchmark method is extremely high. Therefore, the aim of the rest of the methods is to reduce such computational time while obtaining commitment decisions as close as possible to those obtained by the benchmark method. 

\subsection{Single-bus method (SB)} \label{sec:sin}

In this method, all network constraints are removed and the resulting single-bus unit commitment problem is solved. Obviously, this method is the fastest one, but it may provide unit commitment decisions that are very different from the optimal ones. While this method can provide close-to-optimal solutions for power systems in which line congestions hardly ever occur, its use is prohibitive in any other case.

\subsection{Perfect information method (PI)} \label{sec:pin}

Although this method cannot be implemented in practice, its results are useful to understand the discussion in Section~\ref{sec:introduction} on the different types of constraints in MILP problems. This method assumes perfect information on the line capacity constraints that are \textit{active} at the optimal solution and removes all the constraints that are not active, that is, \textit{redundant}, \textit{inactive} and \textit{quasi-active} constraints. In the simulations presented in this paper, the perfect information on the line capacity constraints is obtained from the benchmark method of Section~\ref{sec:ben}. In principle, one may think that this method should provide the same optimal solution as the benchmark's. However, this is not true because of the \textit{quasi-active} constraints defined in Section~\ref{sec:introduction}. Remember that removing \textit{quasi-active} constraint \eqref{ex_c2} from \eqref{ex} provides an optimal solution different from that of the non-reduced problem.

\subsection{Naive method (NV)} \label{sec:nai}

The naive method removes the capacity constraints of those lines that have not been congested for any hour of the available historical data set. In plain words, the intuition of this method is the following: ``If a line has never been congested in the past, it will never get congested in the future''. \textcolor{black}{Despite the low number of constraints that this method is expected to remove, it is a simple method to benchmark the ability of more sophisticated methods to screen out network constraints.} Unlike the three simple methods previously presented, historical data is required to apply this naive method.

\subsection{Constraint generation method (CG)} \label{sec:congen}

Like SB, this method starts by solving the UC without any network constraint. Then, power flows through all lines are computed and the network constraints of those exceeding their capacity limits are iteratively added to the UC problem \cite{Fu2013}. After some iterations, this method provides the same solution as that of BN. The computational burden of this method is highly dependant on the number of iterations required since a reduced UC problem must be solved at each iteration.

\subsection{Zhai's method (ZH)} \label{sec:zha}

This is the method proposed in \cite{Zhai2010}. For each line $l'$, we solve one minimization and one maximization problem with the same objective function and constraints. These problems are jointly formulated in \eqref{zhai}.
\begin{subequations}
\begin{align}
    \min_{p_g,q_{n}}/\max_{p_g,q_{n}}   & \quad \sum_n a_{l'n} q_{n}  \label{zhai_of}\\
    \text{s.t.} & \quad \eqref{uc2_net}, \eqref{uc2_bal}, \eqref{uc2_gen} 
\end{align} \label{zhai}%
\end{subequations}
Optimization problems \eqref{zhai} are different from \eqref{uc2} in the following aspects: i) the objective function aims at maximizing or minimizing the flow through a given line instead of minimizing the total production costs; ii) binary variables are removed and the power output of generators can vary now from 0 to its maximum capacity; iii) line capacity constraints are not considered. Therefore, the feasible region of \eqref{zhai} is larger than that of \eqref{uc2}. Let us denote the optimal solution to the minimization and maximization problems corresponding to line $l'$ as $f^{\min}_{l'}$ and $f^{\max}_{l'}$, respectively. The authors of \cite{Zhai2010} prove that if $|f^{\min}_{l'}|$ and $|f^{\max}_{l'}|$ are lower than  $\overline{f}_{l'}$, then the  capacity constraints of line $l'$ can be removed from \eqref{uc2} without affecting its optimal solution. 

Regarding its applicability, note that optimization problems \eqref{zhai} must be solved for each line whose status we want to infer. For instance, if hourly commitment decisions have to be determined for a system with $N_L$ lines for each hour of the test data set, optimization problem \eqref{zhai} must be solved $2 \times \widehat{N_T} \times N_L$ times. To make it suitable for practical applications, the authors of \cite{Zhai2010} propose an efficient method to find the solution to \eqref{zhai}. 

Albeit not considered in  \cite{Zhai2010}, network constraints \eqref{uc2_cap} can be included in problem \eqref{zhai} too. This way, this method becomes less conservative and removes a higher amount of network constraints. We refer to this variant of the method as ZH+.

\subsection{Roald's method (RO)} \label{sec:roa}

This method has been proposed in \cite{Roald2019}. Similarly to \cite{Zhai2010}, it is based on the solution of one minimization and one maximization problem for each line $l'$. These problems are also jointly formulated as in \eqref{roald}:
\begin{subequations}
\begin{align}
    \min_{p_g,q_{n},d_{n}}/\max_{p_g,q_{n},d_{n}}   & \quad \sum_n a_{l'n} q_{n} \label{roald_of}\\
    \text{s.t.} & \quad \eqref{uc2_net}, \eqref{uc2_bal}, \eqref{uc2_gen}, \eqref{uc2_cap} \\
    & \quad \underline{d}_n \leq d_{n} \leq \overline{d}_n, \quad \forall n \label{roald_dem}
\end{align} \label{roald}%
\end{subequations}
There are two main differences between optimization problems \eqref{zhai} and  \eqref{roald}. First, model \eqref{zhai} considers a fixed demand $d_n$, while  the demand at each bus is a decision variable that can vary from $\underline{d}_n$ to $\overline{d}_n$ in model \eqref{roald}. Second, the capacity factor $\rho_g$ needs to be set in \eqref{roald} to, for example, the maximum capacity factor of each unit $g$. 
If $f^{\min}_{l'}$ and $f^{\max}_{l'}$ correspond to the optimal solutions of \eqref{roald} for a given line $l'$ and it is satisfied that $|f^{\min}_{l'}|<\overline{f}_{l'}$ and $|f^{\max}_{l'}|<\overline{f}_{l'}$, the authors of \cite{Roald2019} propose to remove the capacity constraints of such a line from \eqref{uc2}. 

Although the differences between ZH and RO methods seem small, they involve fundamental consequences regarding their practical applicability. Indeed, let us assume that $\overline{d}_n$ and $\underline{d}_n$ are set to the maximum and minimum demand values at each node during the past year. In that case, having $|f^{\min}_{l'}|<\overline{f}_{l'}$ and $|f^{\max}_{l'}|<\overline{f}_{l'}$ means that line $l'$ is not expected to be congested for a wide range of load variations and its corresponding capacity constraints can be removed for any unit commitment problem to be solved during the next year, for example. This way, if hourly commitment decisions have to be determined for a system of $N_L$ lines for each hour of the test data set, optimization problem \eqref{roald} only has to be solved $2 \times N_L$ times.

\subsection{Data-driven method (DD)} \label{sec:dat}

As discussed, CG, ZH and RO methods are optimization-based and require the solution of extra optimization problems to determine the network constraints to be removed. Besides, ZH and RO only remove \textit{redundant} constraints that do not alter the original feasible region. \textcolor{black}{For this reason, the constraints removed by these methods are not altered by changes in the objective function. However,} this strategy may be too conservative and involve low computational time reductions. 

Conversely, the method proposed in this paper uses historical data on the power system to infer the congestion of the transmission lines via statistical learning methods \cite{Hastie2009}. This shows two main advantages. First, we avoid the need for solving a large number of additional optimization problems. \textcolor{black}{Second, we can efficiently leverage available information on the power system to remove not only \textit{redundant} but also some \textit{inactive} and \textit{quasi-active} constraints. As discussed before, however, while a constraint is redundant irrespective of the objective function of the problem to be solved, a constraint being inactive or quasi-active is contingent on that very same objective function. Therefore, the proposed methodology may be more sensitive to changes in that function.}

Although different learning methods could be applied, we opt for the $K$-nearest neighbors ($K$NN) algorithm because of its simplicity, interpretability, and flexibility to account for physical information about the power system within the learning process. In this paper, this non-parametric (supervised) learning algorithm first proposed in \cite{Fix1951} is used for the classification of congested and not congested lines at time period $t$, based on the information provided by the net demand.

Let us first define the net demand $\tilde{d}_n$ as the difference between the electricity consumption and the renewable generation at each node, i.e., 
\begin{equation}
    \tilde{d}_n = d_n - \sum_{g\in\mathcal{G}^R:b_g=n} \rho_g \overline{p}_g
\end{equation}

Vector $\tilde{\textbf{d}}$ includes the net demand for all nodes, that is, $\tilde{\textbf{d}} = (\tilde{d}_{1}, \ldots, \tilde{d}_{N_B})$. Let us also assume that we have access to historical information on the net demand and the status of a given line $l$ for $N_T$ time periods, represented by the set of pairs $\{(\tilde{\textbf{d}}_t, s_{lt})$, $t = 1, \ldots, N_T\}$. Hence, for each line $l$, and based on the historical data provided by $(\tilde{\textbf{d}}_t, s_{lt})$, $\forall t$, the goal of the $K$NN classification method is to determine the status of line $l$ at a new unseen time period $\hat{t}$, according to the plurality vote of its closest $K$ neighbors. The set of the closest $K$ neighbors is denoted by $\mathcal{N}_K\subset \{1, \ldots, N_T\}$, and is formed by the $K$ most similar time periods. Such a similarity is measured in terms of a distance between the features of $\tilde{\textbf{d}}_t$ and $\tilde{\textbf{d}}_{\hat{t}}$ computed as
\begin{equation}\label{eq: distance definition}
        \text{dist}(\tilde{\textbf{d}}_t,\tilde{\textbf{d}}_{\hat{t}}) = \| \bm{a}_l^T (\tilde{\textbf{d}}_t-\tilde{\textbf{d}}_{\hat{t}}) \|_2
\end{equation}

The distance function \eqref{eq: distance definition} is chosen to be the weighted $\ell_2-$norm due to its well-known performance and simplicity. However, the application of our methodology to other distances such as Minkowski's is straightforward. Notice also that the weights in \eqref{eq: distance definition} are represented by the PTDF of the line under study. Hence, two individuals are close if the net demand of those buses that have a higher impact on the power flow through line $l$ are similar enough.

For each line, $l$, the application of $K$NN to our problem can be thus summarized as follows:
\begin{itemize}
    \item[1)] Use \eqref{eq: distance definition} to compute the distance between the new individual $\hat{t}$ and all existing individuals $t$. 
    \item[2)] Determine the subset $\mathcal{N}_K$ that includes the $K$ time periods with lowest distance $\text{dist}(\tilde{\textbf{d}}_t,\tilde{\textbf{d}}_{\hat{t}})$.
    \item[3)] If $s_{lt}=0\; \forall t \in \mathcal{N}_K$ (which means that line $l$ is not congested for any neighbor), then we assume that the line is not congested for the new individual $\hat{t}$ and its corresponding line capacity constraints \eqref{uc2_cap} can be removed from \eqref{uc2}. Otherwise, such constraints are kept.
\end{itemize}

Note that the number of neighbors $K$ relates somehow with the confidence level of removing the line capacity constraints while obtaining unit commitments that are close enough to those obtained in the original TC-UC problem. The larger $K$, the more conservative we are, and therefore the more accurate the provided solution is. However, increasing $K$ also involves higher computational times. In the limit, if $K = 1$, one line capacity constraint is removed provided that such a line is not congested for the closest neighbor. It is obvious that such a strategy may be quite risky and yield low accurate solutions. On the other hand, if $K$ is equal to the number of individuals, then the proposed method coincides with the naive method of Section \ref{sec:nai} and the only network constraints that can be removed are those that have never been active in the past.

Even though augmenting the number of neighbors increases the conservativeness of the proposed method, there may still be cases where retrieving the original solution is not possible for a sufficiently large value of $K$. For instance, that might be the case of a certain network constraint that only becomes quasi-active or inactive for all historical data points. Regardless of the value of $K$, such network constraint is always removed by the proposed method, which may lead to inaccurate solutions for some operating conditions. If commitment decisions are required to be exactly the same as those obtained by the original TC-UC problem, we propose to equip our data-driven algorithm with a post-processing step based on the constraint generation method as follows: 
\begin{enumerate}
\item  Use the data-driven method to determine the set of network constraints that can be removed.
\item  Solve the reduced unit commitment problem.
\item  Compute the power flow through all lines and determine those that exceed their capacities.
\item  Add the constraints of such lines to the TC-UC problem.
\item  Repeat steps 3 and 4 until all power flows lie between their technical limits.
\end{enumerate}
We refer to this variant of our method as DD+CG. Obviously, the computational time of DD+CG is higher than that of the purely data-driven method DD. Therefore, choosing between these methods depends on the trade-off between accuracy and computational burden.

\section{Comparison} \label{sec:comparison}

The methods presented in Section \ref{sec:methodology} aims at selecting the network constraints that can be removed in order to reduce the computational time required to determine commitment decisions $u_g$ for given values of demand and capacity factors. However, we would also like the commitment decisions $u_g$ yielded by each method to be as close as possible to those obtained by the benchmark method. Therefore, the different methods are compared in terms of the computational time required to solve the reduced TC-UC problem and the accuracy of the obtained solutions. To do so, we proceed as follows:
\begin{itemize}
    \item[1)] Given the historical data, determine the subset of network constraints that can be removed according to each method $m$ described in Section \ref{sec:methodology}, where  $m\in\{\text{BN}, \text{SB}, \text{CG}, \text{PI}, \text{NV}, \text{ZH}, \text{RO}, \text{DD}\}$. The set of removed network constraints is denoted by $\mathcal{L}^m$, and the percentage of removed network constraints by \textcolor{black}{$\mathcal{R}^m = 100\cdot\mathcal{L}^m/N_L$}. The computational burden to determine if the capacity constraints of a given transmission line should be removed (denoted by $\mathcal{T}_1^m$) highly depends on each method. For instance, the time required for methods BN, SB, CG, PI and NV is negligible. Conversely, methods ZH and RO require the resolution of $2 \times \widehat{N_T}$ and two linear programming problems, respectively. Finally, the pre-processing time of the proposed approach based on the $K$NN algorithm is proportional to $N_T \times \widehat{N_T} \times (N_B + K)$.
    \item[2)] Solve the reduced TC-UC problem \eqref{uc2} on the set of unseen time periods $\hat{t}=1,\dots,\widehat{N_T}$ and denote its computational burden as $\mathcal{T}_2^m$.
    \item[3)] Fix commitment decisions to those obtained in Step 2) and solve problem \eqref{uc2} including all network constraints. 
    %
%
    \item[4)] Compute the production cost for step 3) as $\mathcal{C}^m = \sum_g c_gp_{g}$ and the relative infeasibility with respect to the total demand as $\mathcal{I}^m = 100\frac{\sum_n |\epsilon_{n}|}{\sum_n d_{n}}$. This relative infeasibility measures the order of magnitude of the expected deviations between generation and demand during the real-time operation of the system.
\end{itemize}

Methods of Section \ref{sec:methodology} aim at obtaining a value of $\mathcal{T}_1^m + \mathcal{T}_2^m$ as low as possible and a value of $\mathcal{C}^m$ and $\mathcal{I}^m$ as close as possible to those of the BN method.

\section{Illustrative Example} \label{sec:example}

\subsection{Data}

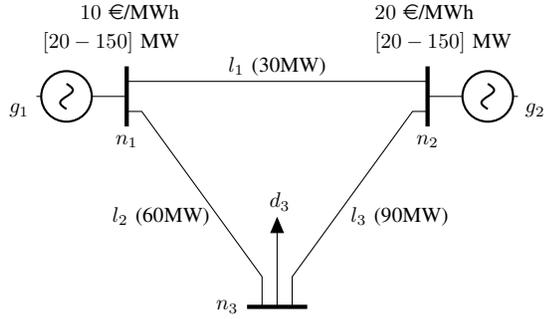
\begin{figure}
\centering
\begin{circuitikz}[scale=0.8,>=triangle 45,every node/.style={scale=0.8}]
\draw [ultra thick] (0,4)  -- (0,3) node[anchor=north]{$n_{1}$};
\draw [ultra thick] (5,4)  -- (5,3) node[anchor=north]{$n_{2}$};
\draw [ultra thick] (2,0) node[anchor=east]{$n_{3}$} -- (3,0);
\draw[] (-1.5,3.5) node[anchor=north east]{$g_1$} to [sV] (-0.5,3.5);
\draw (-0.5,3.5) -- (0,3.5);
\draw[] (6.5,3.5) node[anchor=north west]{$g_2$} to [sV] (5.5,3.5); 
\draw (5.5,3.5) -- (5,3.5); 
\draw[->] (2.5,0) -- +(0,1.5); 
\draw (0,3.75) -- (5,3.75);
\node[] at (2.5,4) {$l_1$ (30MW)};
\draw (0,3.25) -- (0.25,3.25) -- (2.25,0.5) -- (2.25,0); 
\node[anchor=east] at (1.5,1.5) {$l_2$ (60MW)};
\draw (5,3.25) -- (4.75,3.25) -- (2.75,0.5) -- (2.75,0);
\node[anchor=west] at (3.6,1.5) {$l_3$ (90MW)};
\node[anchor=east] at (1,4.9) {$10$ \euro/MWh};
\node[anchor=east] at (1,4.4) {$[20-150]$ MW};
\node[anchor=west] at (4,4.9) {$20$ \euro/MWh};
\node[anchor=west] at (4,4.4) {$[20-150]$ MW};
\node[anchor=west] at (2.25,1.75) {$d_3$};
\end{circuitikz}
\caption{Three-node illustrative example}
\label{fig:3bus}
\end{figure}

To illustrate the differences among the methods of Section~\ref{sec:methodology}, we consider a three-node system with two thermal generating units and one load as depicted in Fig. \ref{fig:3bus}. The minimum and maximum power output of both units are 20MW and 150MW. However, the production cost of units at nodes $n_1$ and $n_2$ is 10 \euro/MWh and 20 \euro/MWh, respectively. The susceptance and capacity of lines $l_1$, $l_2$ and $l_3$ are (1 p.u., 30MW), (2 p.u., 60MW) and (3 p.u., 90MW), in that order. 

Table \ref{tab:data_example} provides the optimal power outputs for six different historical net demand levels together with the status of each line. As observed, line $l_2$ gets congested for demand levels 110, 130 and 150 MW, while $l_3$ only reaches its capacity for the highest value of demand.

\begin{table}[]
    \centering \caption{Illustrative example --- Historical data}
\begin{tabular}{cccccc}
    \hline
    $d_3$ & $p_1$ & $p_2$ & $s_1$ & $s_2$ & $s_3$ \\
    \hline
    50  & 50&0  & 0 & 0 & 0 \\ 
    70  & 70&0  & 0 & 0 & 0 \\
    90  & 70&20 & 0 & 0 & 0 \\
    110 & 73&37 & 0 & 1 & 0 \\
    130 & 66&64 & 0 & 1 & 0 \\
    150 & 60&90 & 0 & 1 & 1 \\
    \hline
    \end{tabular}
    \label{tab:data_example}
\end{table}

\subsection{Results}

Tables \ref{tab:example_85} and \ref{tab:example_125} provide the results after running all the methods described in Section \ref{sec:methodology} for this illustrative example for two demand levels of 85 MW and 125 MW, respectively. Columns 2 and 3 indicate the removed network constraints ($\mathcal{L}^m$) and the set of committed units for each method $m$. For such commitment decisions, columns 4-5, 6-8 and 9 include the dispatch quantities, the flows through the lines and the slack variable at bus 3 if all network constraints are considered. We discuss below the performance of each method:


\begin{table}[]
    \centering
    \setlength{\tabcolsep}{4pt}
    \caption{Illustrative example --- Comparison of methods ($d_3 = 85$MW)}
    \begin{tabular}{l|cc|cccccc}
    \hline
    $m$ & $\mathcal{L}^m$ & $u_g=1$ & $p_{g_1}$ & $p_{g_2}$ & $f_{l_1}$ & $f_{l_2}$ & $f_{l_3}$ & $\epsilon_3$ \\
    \hline
    \textbf{BN} & $-$ & $g_1,g_2$ & 65 & 20 & 14.1 & 50.9 & 34.1 & 0 \\
    SB & $l_1,l_2,l_3$ & $g_1$ & 82.5 & 0 & 22.5 & 60.0 & 22.5 & -2.5  \\
    PI & $l_1,l_2,l_3$ & $g_1$ & 82.5 & 0 & 22.5 & 60.0 & 22.5 & -2.5  \\
    \textbf{NV} & $l_1$ & $g_1,g_2$ & 65 & 20 & 14.1 & 50.9 & 34.1 & 0  \\
    \textbf{CG} & $l_1,l_3$ & $g_1,g_2$ & 65 & 20 & 14.1 & 50.9 & 34.1 & 0 \\
    \textbf{ZH} & $l_1,l_3$ & $g_1,g_2$ & 65 & 20 & 14.1 & 50.9 & 34.1 & 0 \\
    \textbf{RO} & $-$ & $g_1,g_2$ & 65 & 20 & 14.1 & 50.9 & 34.1 & 0 \\
    DD2 & $l_1,l_2,l_3$ & $g_1$ & 82.5 & 0 & 22.5 & 60.0 & 22.5 & -2.5 \\
    \textbf{DD3} & $l_1,l_3$ & $g_1,g_2$ & 65 & 20 & 14.1 & 50.9 & 34.1 & 0 \\
    \textbf{DD6} & $l_1$ & $g_1,g_2$ & 65 & 20 & 14.1 & 50.9 & 34.1 & 0 \\
    \hline
    \end{tabular}
    \label{tab:example_85}
\end{table}

\begin{table}[]
    \centering
    \setlength{\tabcolsep}{4pt}
    \caption{Illustrative example --- Comparison of methods ($d_3 = 125$MW)}
    \begin{tabular}{l|cc|cccccc}
    \hline
    $m$ & $\mathcal{L}^m$ & $u_g=1$ & $p_{g_1}$ & $p_{g_2}$ & $f_{l_1}$ & $f_{l_2}$ & $f_{l_3}$ & $\epsilon_3$ \\
    \hline
    \textbf{BN }& $-$ & $g_1,g_2$ & 68 & 57 & 8.2 & 60.0 & 65.2 & 0 \\
    SB & $l_1,l_2,l_3$ & $g_1$ & 82.5 & 0 & 22.5 & 60.0 & 22.5 & -42.5 \\
    PI & $l_1,l_3$ & $g_1,g_2$ & 68 & 57 & 8.2 & 60.0 & 65.2 & 0 \\
    \textbf{NV} & $l_1$ & $g_1,g_2$ & 68 & 57 & 8.2 & 60.0 & 65.2 & 0 \\
    \textbf{CG }& $l_1,l_3$ & $g_1,g_2$ & 68 & 57 & 8.2 & 60.0 & 65.2 & 0 \\
    \textbf{ZH} & $-$ & $g_1,g_2$ & 68 & 57 & 8.2 & 60.0 & 65.2 & 0 \\
    \textbf{RO} & $-$ & $g_1,g_2$ & 68 & 57 & 8.2 & 60.0 & 65.2 & 0  \\
    DD2 & $l_1,l_3$ & $g_1,g_2$ & 68 & 57 & 8.2 & 60.0 & 65.2 & 0 \\
    \textbf{DD3} & $l_1$ & $g_1,g_2$ & 68 & 57 & 8.2 & 60.0 & 65.2 & 0 \\
    \textbf{DD6} & $l_1$ & $g_1,g_2$ & 68 & 57 & 8.2 & 60.0 & 65.2 & 0 \\
    \hline
    \end{tabular}
    \label{tab:example_125}
\end{table}

\begin{itemize}
    \item Benchmark method (BN): This method considers all constraints and involves the highest computational burden for both demand levels.
    \item Single-bus method (SB): This method disregards all network constraints and thus only the cheapest unit $g_1$ is committed for both demand levels. For such commitment decisions, the slack variable at bus 3 amounts to -2.5MW and -42.5MW to keep all power flows within their limits for $d_3=85$MW and $d_3=125$MW, respectively.
    \item Perfect information method (PI): Particularly interesting are the results corresponding to the PI method for $d_3=85$MW. Although none of the line capacity constraints are binding for this demand level in the BN method, if all of them are removed, the obtained commitment decisions are different than those of BN. This occurs because the capacity constraint of line $l_2$ is, for $d_3=85$MW, a \textit{quasi-active} constraint, as defined in Section \ref{sec:introduction}. For $d_3=125$MW, this method keeps the capacity constraint of $l_2$ and obtains the same solution as BN.
    \item Naive method (NV): Since $l_1$ is the only line that has not been congested in the past, the NV method keeps the constraints of the other two lines for both demand levels and yields the same solution as that of BN. 
    \item Constraint generation method (CG): For this example, this method requires two iterations. In the first one, all network constraints are disregarded and in the second one, the network constraints corresponding to $l_2$ are included. All results coincide with that of BN.
    \item Zhai's method (ZH): If optimization problem \eqref{zhai} is solved for $d_3=85$MW, only the flow through line $l_2$ exceeds its capacity. For $d_3=125$MW, the three lines reach their capacities and thus cannot be removed from \eqref{uc2}. For both demands, ZH provides the same solution as that of BN.
    \item Roald's method (RO): If optimization problem \eqref{roald} is solved for each line with $d_3$ varying between its minimum and maximum historical values ($50\leq d_3 \leq 150$), all of them exceed their capacities and therefore, no constraint is removed from \eqref{uc2}. The results corresponding to RO methods coincide then with those of BN. 
    \item Data-driven method (DD): We include the results of the proposed data-driven method for 2, 3 and 6 neighbors and denote them as DD2, DD3, and DD6, respectively. As this example only includes one feature (electricity demand at bus 3), considering weights in the proposed method does not apply. For $d_3=85$MW, the two closest neighbors are the data points corresponding to $d_3=70$MW and $d_3=90$MW. Since the three lines are uncongested for both neighbors, the DD2 method does not consider any network constraints. Since the capacity constraint of line $l_2$ for this demand level is a \textit{quasi-active} constraint, removing it from \eqref{uc2} derives commitment decisions different from those of BN. For the same demand level, the three closest neighbors correspond to $d_3=70$MW, $d_3=90$MW and $d_3=110$MW. Notice that $l_2$ is congested for $d_3=110$MW and hence, DD3 does not remove the capacity constraint of such a line and the obtained commitment coincides with that of BN. This illustrates that increasing the number of neighbors $K$ in the proposed data-driven method reduces the probability of removing \textit{quasi-active} constraints in the reduced TC-UC problem. For $d_3=125$MW, the two closest neighbors correspond to $d_3=110$MW and $d_3=130$MW, which implies that DD2 keeps the capacity constraint of $l_2$. Similarly, the three closest neighbors are $d_3=110$MW $d_3=130$MW and $d_3=150$MW and thus, DD3 only removes the capacity constraint of $l_1$. For DD2 and DD3, the obtained commitment coincides with that of BN. Finally, if the number of neighbors is equal to the total number of historical data, the results provided by DD6 coincide with those given by NV for both demand levels. 
\end{itemize}

All the methods that provide the same commitment decisions as those of the BN method for both demand levels are highlighted in bold font in Tables \ref{tab:example_85} and \ref{tab:example_125}. These methods are compared next in terms of its ability to reduce the computational burden of the original TC-UC problem. Given the reduced size of this example, comparing computational times is not conclusive and therefore, we assume that the number of line capacity constraints removed for both demand levels is directly related to the time reduction achieved by each method. In that sense, RO involves the same computational time as BN since no constraints are removed. Method CG removes four constraints out of six, but it has to solve the UC twice (once without any constraint, and once with the constraints corresponding to lines $l_1$ and $l_3$), which significantly increases the computational time. Methods NV, ZH, and DD6 relatively reduce the computational burden since two out of six constraints (three constraints per each demand level) are removed. Finally, the proposed method DD3  achieves the highest computational time reduction since three out of six constraints are removed.

\section{Simulation Results} \label{sec:num_results}

This section provides simulation results from a medium- and a large-size power system with 73 and 2000 buses, respectively. The simulations have been performed on a Linux-based
server with one CPU clocking at 2.6 GHz and 2 GB of RAM using CPLEX 12.6.3 under Pyomo 5.2. All provided times are CPU times.

\subsection{Medium-Size Case Study} \label{sec:case_study}

In this section, we analyze and discuss the results of all the methods presented in Section \ref{sec:methodology} for the IEEE RTS-96 test system \cite{Grigg1999a} modified to accommodate 19 wind farms as proposed in \cite{Pandzic2015}. This power system includes 73 nodes and 120 transmission lines. Technical information on generating units and transmission lines is available in \cite{IEEE96_2019}. This repository also includes data corresponding to hourly electricity demand and wind generation for 360 days. Finally, the hourly status of each line (1 if congested and 0 otherwise) obtained after solving the TC-UC problem \eqref{uc2} are also provided.

The whole data set is split into two subsets. The \textit{training data set} includes the first 300 days of the year, while the \textit{test data set} comprises the last 60 days. Different uses are given to the training subset depending on the method applied. For instance, the naive method removes the capacity constraints of the lines that have never become congested in this data subset. Roald's method uses this information to determine the values of $\underline{d}_n$, $\overline{d}_n$ and $\rho_g$. Finally, the data-driven method proposed in this paper chooses the closest neighbors of a new data point from this data subset. On the other hand, the test data subset is used to compare the performance of the different methods as explained in Section \ref{sec:comparison}. 

To assess the performance of the proposed data-driven method, we have carried out numerical simulations for three different levels of network congestion, namely, a \textit{medium-congested case} (MC) with the original line capacities, a \textit{low-congested case} (LC) obtained by multiplying the line capacities by two, and a \textit{high-congested case} (HC) obtained by dividing the line capacities by two. The number of lines that get ever congested for the whole year for LC, MC, and HC are 1, 8, and 39, respectively. Besides, the line that most often gets congested for LC, MC and HC reaches its corresponding capacity limit during 316, 3714 and 5727 hours of the year. The time required by the BN method to solve the TC-UC problem for the 60 days of the test subset for the LC, MC, and HC cases amounts to 389, 235 and 1153 seconds, in that order. This means that reducing the computational burden becomes more imperative for highly congested power systems. We next analyze the results for each case separately.

\subsubsection{Medium-congested case}

For all the methods compared in this paper, Table \ref{tab:case_study_med} summarizes the main results, namely, the percentage of line capacity constraints that are removed by each method ($\mathcal{R}^m$), the relative cost error with respect to the benchmark method ($\Delta \mathcal{C}^m = 100(\mathcal{C}^m-\mathcal{C}^{\text{BN}})/\mathcal{C}^{\text{BN}}$), the relative infeasibility of each method ($\mathcal{I}^m$), computational times $\mathcal{T}_1^m,\mathcal{T}_2^m$, and the relative computational time with respect to the benchmark method $ \left( \tau^m = 100 \cdot (\mathcal{T}_1^m + \mathcal{T}_2^m) / \mathcal{T}_2^{\text{BN}} \right)$ . For method RO100, the upper and lower bounds on the variable $d_n$ are set to the maximum and minimum historical electricity demand at each node. The capacity factor $\rho_g$ is fixed to the maximum value as well. For RO95 and RO90, the demand bounds are set to percentiles 95\% and 5\%, and 90\% and 10\%, respectively. Similarly, $\rho_g$ is set to the percentile 95\% and 90\% of the historical data. It is also worth mentioning that the results provided in Table \ref{tab:case_study_med} are accumulated values for the 60 days of the test data set. The optimality gap is set to 0\% for this medium-size case study.

\begin{table}[]
    \setlength{\tabcolsep}{4pt}
    \centering \caption{Medium-congested 73-bus system --- Comparison of methods}
    \begin{tabular}{lcccccc}
    \hline
    Method($m$) & $\mathcal{R}^m$(\%) & $\Delta \mathcal{C}^m$(\%) & $\mathcal{I}^m$(\%) & $\mathcal{T}_1^m$(s) & $\mathcal{T}_2^m$(s) & $\tau^m$(\%) \\
    \hline
    BN & 0 & 0.00 & 0.000 & 0.0 & 388.6 & 100.0 \\
    SB & 100 & -4.08 & 3.390 & 0.0 & 90.1 & 23.1 \\
    PI & 99.3 & -0.11 & 0.145 & 0.0 & 95.5 & 24.6 \\
    NV & 94.2 & 0.00 & 0.000 & 0.0 & 131.8 & 33.9 \\
    CG & 98.7 & 0.00 & 0.000 & 0.0 & 273.7 & 70.4 \\
    ZH & 64.3 & 0.00 & 0.000 & 442.2 & 119.8 & 144.6 \\
    ZH+ & 82.2 & 0.00 & 0.000 & 647.3 & 76.2 & 186.2 \\
    RO100 & 53.3 & 0.00 & 0.000 & 2.1 & 246.8 & 64.1 \\
    RO95 & 63.3 & 0.00 & 0.000 & 2.1 & 197.7 & 51.4 \\
    RO90 & 67.5 & 0.00 & 0.000 & 2.1 & 198.0 & 51.5 \\
    DD5 & 99.2 & 0.00 & 0.002 & 0.4 & 110.5 & 28.6 \\
    DD50 & 98.9 & 0.00 & 0.000 & 0.4 & 112.7 & 29.1 \\
    DD500 & 98.3 & 0.00 & 0.000 & 0.4 & 99.1 & 25.6 \\
    DD5+CG & 99.2 & 0.00 & 0.000 & 0.5 & 136.0 & 35.1 \\
    DD50+CG & 98.8 & 0.00 & 0.000 & 0.4 & 114.6 & 29.6 \\
    DD500+CG & 98.2 & 0.00 & 0.000 & 0.5 & 100.2 & 25.9 \\
    \hline
    \end{tabular}
    \label{tab:case_study_med}
\end{table}

The BN method does not remove any network constraints and as a result, obtains the optimal solution requiring a high computational time. The SB method removes all network constraints, which significantly reduces the computational time but provides commitment solutions that are quite different from those corresponding to BN as shown by the values of $\Delta \mathcal{C}^m$ and $\mathcal{I}^m$. Interestingly, although the PI method only removes the capacity constraints of the lines that are not congested, it also provides inaccurate commitment decisions due to the \textit{quasi-active} constraints defined in Section \ref{sec:introduction}. The NV method removes 94.2\% of the network constraints and reduces the computational time to 33.9\% of the time required by the BN method. The CG method only adds 1.3\% of the network constraints, but the time reduction is not very significant since the UC has to be solved several times.

As previously discussed, the ZH method is quite conservative and as such, only removes 64.3\% of the line capacity constraints, which leads to a relative computational time of 144.6\%. That is, if the time required to solve optimization problems \eqref{zhai} is considered, the total computational time exceeds that of the original TC-UC problem. Although method ZH+ is less conservative and removes a higher amount of network constraints, it does not involve any computational savings either. RO method also removes a low percentage of capacity constraints, although such a number increases as bounds on $d_{n}$ are tightened. This method reduces the computational burden by half, approximately. Finally, the proposed data-driven method removes a number of constraints above 98\% for values of $K$ equal to 5, 50 and 500 and reduces the computational time to 25-30\% of the time required by the BN method. Notice that DD5 involves unit commitment decisions slightly different from those of the BN method since the relative infeasibility amounts to 0.002\%. As stated at the end of Section \ref{sec:dat}, setting $K$ to a larger value increases the chances of obtaining more accurate solutions. In effect, if $K$ is increased to 50 or 500, the proposed method provides the same unit commitment decisions as BN at the lowest computational time. Alternatively, such inaccuracies can also be removed by combining DD5 with CG at the expense of slightly increasing the computational time. Here, the potential user faces a clear trade-off between the computational burden and the accuracy of the obtained solutions.

\subsubsection{Low-congested case}

Table \ref{tab:case_study_low} provides results similar to those in Table \ref{tab:case_study_med} for the low-congested case. First, it can be observed that the results given by the SB method are not as inaccurate as in the medium-congested case. Indeed, since line congestion rarely ever happens, removing all network constraints is a sensible strategy that involves an 80\% time reduction. However, similar time reductions are also achieved by the other more sophisticated methods such as RO or DD. Therefore, for low-congested systems, simple methods such as SB or NV perform good enough to reduce the computational complexity of the TC-UC problem without significantly affecting the accuracy of the unit commitment decisions.

\begin{table}[]
    \setlength{\tabcolsep}{4pt}
    \centering \caption{Low-congested 73-bus system ---  Comparison of methods}
    \begin{tabular}{lcccccc}
    \hline
    Method($m$) & $\mathcal{R}^m$(\%) & $\Delta \mathcal{C}^m$(\%) & $\mathcal{I}^m$(\%) & $\mathcal{T}_1^m$(s) & $\mathcal{T}_2^m$(s) & $\tau^m$(\%) \\
    \hline
    BN & 0 & 0.00 & 0.000 & 0.0 & 234.9 & 100 \\
    SB & 100 & -0.04 & 0.044 & 0.0 & 46.6 & 19.8 \\
    PI & 100 & 0.01 & 0.001 & 0.0 & 37.2 & 15.9 \\
    NV & 99.2 & 0.00 & 0.000 & 0.0 & 38.7 & 16.5 \\
    CG & 100 & 0.00 & 0.000 & 0.0 & 63.5 & 27.0 \\
    ZH & 94.8 & 0.00 & 0.000 & 341.8 & 25.1 & 156.2 \\
    ZH+ & 95.6 & 0.00 & 0.000 & 377.2 & 22.8 & 170.3 \\
    RO100 & 85.8 & 0.00 & 0.000 & 0.8 & 57.5 & 24.8 \\
    RO95 & 90.0 & 0.00 & 0.000 & 0.7 & 51.0 & 22.0 \\
    RO90 & 91.7 & 0.00 & 0.000 & 0.8 & 46.7 & 22.2 \\
    DD5 & 99.9 & 0.01 & 0.002 & 0.0 & 44.0 & 18.7 \\
    DD50 & 99.8 & 0.00 & 0.000 & 0.0 & 45.3 & 17.5 \\
    DD500 & 99.5 & 0.00 & 0.000 & 0.0 & 48.5 & 19.3 \\
    DD5+CG & 99.9 & 0.00 & 0.000 & 0.0 & 48.5 & 20.7 \\
    DD50+CG & 99.8 & 0.00 & 0.000 & 0.0 & 53.0 & 22.6 \\
DD500+CG & 99.5 & 0.00 & 0.000 & 0.0 & 51.6 & 22.0 \\
    \hline
    \end{tabular}
    \label{tab:case_study_low}
\end{table}

\subsubsection{High-congested case}

For the high-congested case, Table \ref{tab:case_study_hig} shows that the results of the compared methods are more different in terms of computational burden and accuracy. For instance, the SB method provides TC-UC decisions that are quite inaccurate since the production cost and relative infeasibility are very different from those of BN. Notice that the NV method provides acceptable results, with a relative computational time of 33.3\%, almost no infeasibilities and a small cost error. The CG method provides the same results as BN but the computational savings are not very notable due to the high number of iterations required. For the high-congested case, methods ZH and ZH+ are also outperformed by other methods such as CG, RO, DD, and DD+CG. Due to the high level of congestion, method RO removes a very low number of network constraints and therefore, leads to quite small computational savings. The proposed DD method is able to efficiently learn from available data and remove between 83.9\% and 93.1\% of the line capacity constraints. This results in time reductions that are between 75\% and 90\% for the three values of $K$ here considered. The largest computational savings are obtained for DD5 but at the expense of getting significant values of relative infeasibility and cost error. On the other hand, the computational time of DD500 is also noticeably low, while keeping $\Delta\mathcal{C}^{m}$ and $\mathcal{I}^{m}$ to values very close to 0\%, which means that the accuracy of the obtained commitment decisions is very high. Therefore, the proposed methodology is able to efficiently reduce the computational burden of the TC-UC problem by a factor of five and provide commitment decisions almost identical to those obtained by the benchmark model. Finally, the inaccuracies are equal to 0 for the three instances of method DD+CG at the expense of reducing the time savings. For example, the method DD50+CG reduces the computational time to 32.7\% while obtaining the same commitment decisions as those of the benchmark.

\begin{table}[]
    \setlength{\tabcolsep}{4pt}
    \centering \caption{High-congested 73-bus system --- Comparison of methods}
    \begin{tabular}{lcccccc}
    \hline
    Method($m$) & $\mathcal{R}^m$(\%) & $\Delta \mathcal{C}^m$(\%) & $\mathcal{I}^m$(\%) & $\mathcal{T}_1^m$(s) & $\mathcal{T}_2^m$(s) & $\tau^m$(\%) \\
    \hline
    BN & 0 & 0.00 & 0.000 & 0.0 & 1153.3 & 100.0 \\
    SB & 100 & -20.14 & 10.557 & 0.0 & 35.5 & 3.1 \\
    PI & 95.1 & 1.65 & 0.435 & 0.0 & 65.7 & 5.7 \\
    NV & 72.5 & 0.01 & 0.001 & 0.1 & 383.7 & 33.3 \\
    CG & 83.2 & 0.00 & 0.000 & 0.0 & 753.9 & 65.4 \\
    ZH & 11.3 & 0.00 & 0.000 & 357.6 & 1100.6 & 126.47 \\
    ZH+ & 49.7 & 0.00 & 0.000 & 409.3 & 587.4 & 86.44 \\
    RO100 & 21.7 & 0.00 & 0.000 & 0.8 & 926.2 & 80.4 \\
    RO95 & 23.3 & 0.00 & 0.000 & 0.8 & 880.7 & 76.5 \\
    RO90 & 25.0 & 0.00 & 0.000 & 0.8 & 900.0 & 78.1 \\
    DD5 & 93.1 & 1.43 & 0.204 & 0.9 & 118.0 & 10.3 \\
    DD50 & 89.4 & 0.79 & 0.107 & 0.9 & 144.1 & 12.6 \\
    DD500 & 83.9 & 0.01 & 0.031 & 0.9 & 254.2 & 22.1 \\
    DD5+CG & 91.7 & 0.00 & 0.000 & 0.9 & 407.6 & 35.4 \\
    DD50+CG & 88.8 & 0.00 & 0.000 & 0.9 & 375.6 & 32.7 \\
    DD500+CG & 83.8 & 0.00 & 0.000 & 0.9 & 535.4 & 46.5 \\
    \hline
    \end{tabular}
    \label{tab:case_study_hig}
\end{table}

\subsection{Large-Size Case Study} \label{sec:case_study_large}

In order to show the efficiency of the proposed approach for realistic cases, this section compares the simulation results of the different methods for a power system in Texas with 2000 buses and 3206 lines \cite{Birchfield2017}. All system data are available at \cite{pglib}. As proposed in \cite{Roald2019}, the minimum output of generating unis is set to the maximum between the value specified in the dataset and 10\% of their generation capacity to obtain more challenging TC-UC problems.

To generate different system operating conditions throughout a whole year, the electricity demand at each bus is randomly sampled from uniform distributions between 0 and twice the nominal demand at each bus. Under such conditions, 10\% of the lines become congested during the year, and the line that most often gets congested reaches its capacity limit during 4000 hours. As in the previous case study, the first 300 days of the year are used as the training data set, and the following 60 days as the test data set. To keep computational times within reasonable limits, optimality gap is set to 1\%.

Table \ref{tab:case_study_large} provides the simulation results. Given the poor performance of ZH and ZH+ for the previous medium-size case study, these two methods have not been considered here.

\begin{table}[]
    \setlength{\tabcolsep}{4pt}
    \centering \caption{2000-bus system --- Comparison of methods}
    \begin{tabular}{lcccccc}
    \hline
    Method($m$) & $\mathcal{R}^m$(\%) & $\Delta \mathcal{C}^m$(\%) & $\mathcal{I}^m$(\%) & $\mathcal{T}_1^m$(s) & $\mathcal{T}_2^m$(s) & $\tau^m$(\%) \\
    \hline
    BN & 0 & 0.00 & 0.052 & 0.0 & 53478.3 & 100.0 \\
    SB & 100 & -2.17 & 0.263 & 0.0 & 235.2 & 0.4 \\
    PI & 99.7 & -0.22 & 0.431 & 0.0 & 546.8 & 1.0 \\
    NV & 92.3 & 0.00 & 0.052 & 0.1 & 5687.4 & 10.6 \\
    CG & 98.8 & 0.00 & 0.052 & 0.0 & 4757.8 & 8.9 \\
    RO100 & 54.3 & 0.00 & 0.052 & 42.7 & 34581.8 & 64.7 \\
    RO95 & 54.9 & 0.00 & 0.052 & 42.9 & 31209.3 & 58.4 \\
    RO90 & 54.9 & 0.00 & 0.052 & 42.6 & 29552.8 & 55.3 \\
    DD5 & 99.5 & -0.08 & 0.309 & 82.0 & 627.9 & 1.3 \\
    DD50 & 98.6 & 0.04 & 0.084 & 89.4 & 1161.3 & 2.3 \\
    DD500 & 96.8 & 0.02 & 0.055 & 106.0 & 2333.2 & 4.6 \\
    DD5+CG & 99.2 & 0.00 & 0.052 & 83.9 & 2720.9 & 5.2 \\
    DD50+CG & 98.5 & 0.00 & 0.052 & 90.1 & 2760.4 & 5.3 \\
    DD500+CG & 96.8 & 0.00 & 0.052 & 106.2 & 3517.4 & 6.8 \\
    \hline
    \end{tabular}
    \label{tab:case_study_large}
\end{table}

Notice that the benchmark method (BN) involves small infeasibilities caused by the insufficient network capacity to satisfy the demand in some hours of the year. \textcolor{black}{This means that there exists no unit commitment which avoids infeasibilities for the 60 days of the test data set.} Methods SB, PI, DD5 and DD50 significantly reduce the computational time, but imply variations of the cost and/or the relative infeasibility with respect to BN. Methods RO100, RO95 and RO90 obtain the same solution as that of BN but only save less than a half of the computational time. Methods NV and CG also obtain the same solution as that of BN and reduce the computational burden to 10.6\% and 8.9\%, thus outperforming RO method. Our purely data-driven method DD500 further reduces the computational time obtained by the RO approach to just a 4.6\% of the time used in the BN method at the expense of yielding commitment decisions that are slightly different from those derived by BN. Therefore, the proposed methodology is able to obtain highly accurate solutions in a significantly reduced amount of time. Finally, adding a post-processing step to the proposed method slightly increases the computational time with respect to the purely data-driven method but achieves solutions equal to that of BN. In fact, the relative computational times of DD5+CG, DD50+CG and DD500+CG amount to 5.2\%, 5.3\% and 6.8\%, which are lower than any other method providing the same commitment decisions as the benchmark's.

\subsection{Discussion}

To properly evaluate the practical implications of the solution provided by the proposed method, it is important to keep in mind that the TC-UC problem is usually not solved to global optimality for real-size power systems and mip gap values between 0.5\% and 1\% are usually accepted \cite{Streiffert2005}. Therefore, the suboptimality of the unit commitment determined by DD is significantly below typical values. Secondly, actual power systems include a series of fast-responding ancillary services to balance small deviations between generation and consumption of electricity caused, for example, by forecast errors. These ancillary services include, for example, spinning reserves and frequency control equipments. Therefore, the small inaccuracies yielded by the proposed method in some instances will be handled by such ancillary services without compromising the stability of the power system operation \cite{Banakar2008}. 

In any case, if retrieving the same commitment decisions as those of the original TC-UC problem is a priority, this can be pursued in two different ways. One possibility consists in increasing the number of neighbors, $K$, in our data-driven proposal. As proven by the simulations, the larger the value of $K$, the lower the inaccuracies of the obtained solutions at the expense of reducing the number of removed network constraints and the computational savings. Alternatively, the proposed method can be combined with the CG method to guarantee that the obtained solution is always equal to that of the BN method as previously discussed. In such a case, the computational burden also increases, albeit slightly in general, if compared to the purely data-driven method.

In summary, methods CG, ZH and RO ensure that the obtained commitment decisions are equal to those of the original TC-UC problem but involve limited computational savings. On the other hand, the proposed method DD significantly reduces the computational burden while yielding commitment decisions that may be a bit different from the optimal ones depending on the network congestion level and the number of neighbors. Finally, method DD+CG is able to retrieve the optimal commitment at the lowest computational time.

\section{Conclusions} \label{sec:conclusions}

Reducing the computational burden of the TC-UC problem is a very relevant research topic within the PES community. One approach to address this issue consists in disregarding the capacity constraints of those transmission lines that will not become active during operation. Although there are some existing works in this direction, most of them focus on removing \textit{redundant} network constraints in order not to alter the original feasible region of the TC-UC problem.

In this paper, we propose a data-driven approach that takes advantage of historical information to disregard both \emph{redundant} and \emph{inactive} network constraints. In doing so, we remove a larger number of network constraints than existing methods and achieve higher computational savings, which range from 70\% to 98\% depending on the congestion level of the power system and its size. As a counterpart, the proposed approach yields commitment decisions that may be slightly different from those obtained by the original TC-UC problem. If required, such small inaccuracies can be completely eliminated by running a post-processing step based on a constraint generation procedure to recover the solution to the original TC-UC problem.

Future research is needed to extend the proposed methodology to the multi-period unit commitment problem, where more detailed technical and economic characteristics of generating units are considered. \textcolor{black}{It would be also interesting to investigate the performance of the proposed approach in the security-constrained unit commitment problem.} Besides, the proposed method can be also extended to screen out other technical constraints such as ramping or capacity limits of generating units. Finally, investigating the performance of different statistical learning techniques to reduce the computational burden of the TC-UC problem also requires further research.


%





\ifCLASSOPTIONcaptionsoff
  \newpage
\fi

\balance



\bibliographystyle{IEEEtran}
\bibliography{references}
\end{document}